\documentclass[11pt]{amsart}
\usepackage{amsmath,amssymb}
\newtheorem{theorem}{Theorem}[section]
\newtheorem{proposition}[theorem]{Proposition}
\newtheorem{corollary}[theorem]{Corollary}
\newtheorem{lemma}[theorem]{Lemma}

\begin{document}

\title[Sobolev inequalities]{Sobolev inequalities in manifolds with nonnegative curvature}
\author{Simon Brendle}
\address{Department of Mathematics \\ Columbia University \\ New York NY 10027}
\begin{abstract}
We prove a sharp Sobolev inequality on manifolds with nonnegative Ricci curvature. Moreover, we prove a Michael-Simon inequality for submanifolds in manifolds with nonnegative sectional curvature. Both inequalities depend on the asymptotic volume ratio of the ambient manifold. 
\end{abstract}
\thanks{This project was supported by the National Science Foundation under grant DMS-1806190 and by the Simons Foundation.}

\maketitle 

\section{Introduction} 

Let $M$ be a complete noncompact manifold of dimension $k$ with nonnegative Ricci curvature. The asymptotic volume ratio of $M$ is defined as 
\[\theta := \lim_{r \to \infty} \frac{|\{p \in M: d(p,q)  < r\}|}{|B^k| \, r^k},\] 
where $q$ is some fixed point in $M$ and $B^k$ denotes the unit ball in $\mathbb{R}^k$. The Bishop-Gromov relative volume comparison theorem implies that the limit exists, and that $\theta \leq 1$. Note that $\theta$ does not depend on the choice of the point $q$.

Our first result gives a sharp Sobolev inequality on manifolds with nonnegative Ricci curvature.

\begin{theorem}
\label{sobolev.inequality.a}
Let $M$ be a complete noncompact manifold of dimension $n$ with nonnegative Ricci curvature. Let $D$ be a compact domain in $M$ with boundary $\partial D$, and let $f$ be a positive smooth function on $D$. Then 
\[\int_D |\nabla f| + \int_{\partial D} f \geq n \, |B^n|^{\frac{1}{n}} \, \theta^{\frac{1}{n}} \, \Big ( \int_D f^{\frac{n}{n-1}} \Big )^{\frac{n-1}{n}},\] 
where $\theta$ denotes the asymptotic volume ratio of $M$.
\end{theorem}

Moreover, we are able to characterize the case of equality in Theorem \ref{sobolev.inequality.a}: 

\begin{theorem}
\label{rigidity.theorem.a}
Let $M$ be a complete noncompact manifold of dimension $n$ with nonnegative Ricci curvature. Let $D$ be a compact domain in $M$ with boundary $\partial D$, and let $f$ be a positive smooth function on $D$. Suppose that 
\[\int_D |\nabla f| + \int_{\partial D} f = n \, |B^n|^{\frac{1}{n}} \, \theta^{\frac{1}{n}} \, \Big ( \int_D f^{\frac{n}{n-1}} \Big )^{\frac{n-1}{n}} > 0,\] 
where $\theta$ denotes the asymptotic volume ratio of $M$. Then $f$ is constant, $M$ is isometric to Euclidean space, and $D$ is a round ball.
\end{theorem}

Putting $f=1$ in Theorem \ref{sobolev.inequality.a}, we obtain a sharp isoperimetric inequality.

\begin{corollary}
\label{isop.inequality.a}
Let $M$ be a complete noncompact manifold of dimension $n$ with nonnegative Ricci curvature. Let $D$ be a compact domain in $M$ with boundary $\partial D$. Then 
\[|\partial D| \geq n \, |B^n|^{\frac{1}{n}} \, \theta^{\frac{1}{n}} \, |D|^{\frac{n-1}{n}},\] 
where $\theta$ denotes the asymptotic volume ratio of $M$.
\end{corollary}

Corollary \ref{isop.inequality.a} is similar in spirit to the L\'evy-Gromov inequality for manifolds with Ricci curvature at least $n-1$ (cf. \cite{Gromov}, Appendix C). The L\'evy-Gromov inequality was recently generalized in \cite{Klartag} and \cite{Cavalletti-Mondino}. 

In the three-dimensional case, Corollary \ref{isop.inequality.a} was proved in a recent work of V.~Agostiniani, M.~Fogagnolo, and L.~Mazzieri (cf. \cite{Agostiniani-Fogagnolo-Mazzieri}, Theorem 6.1). The proof of Theorem 6.1 in \cite{Agostiniani-Fogagnolo-Mazzieri} builds on an argument due to G.~Huisken \cite{Huisken} and uses mean curvature flow. 

We will present the proof of Theorem \ref{sobolev.inequality.a} in Section \ref{proof.of.sobolev.inequality.a}. The proof of Theorem \ref{sobolev.inequality.a} uses the Alexandrov-Bakelman-Pucci method and is inspired in part by an elegant argument due to X.~Cabr\'e \cite{Cabre2} (see also \cite{Cabre1},\cite{Cabre-Ros-Oton-Serra},\cite{Trudinger},\cite{Wang-Zhang},\cite{Xia-Zhang} for related work). The proof of Theorem \ref{rigidity.theorem.a} will be discussed in Section \ref{proof.of.rigidity.theorem.a}. 

%
%

We next turn to Sobolev inequalities for submanifolds. In a recent paper \cite{Brendle}, we proved a Michael-Simon-type inequality for submanifolds in Euclidean space. While the classical Michael-Simon inequality (cf. \cite{Allard}, \cite{Michael-Simon}) is not sharp, our inequality is sharp if the codimension is at most $2$. In particular, the results in \cite{Brendle} imply a sharp isoperimetric inequality for minimal submanifolds in Euclidean space of codimension at most $2$, answering a question first studied by Torsten Carleman \cite{Carleman} in 1921. 

The following theorem generalizes the main result in \cite{Brendle} to the Riemannian setting.

\begin{theorem}
\label{sobolev.inequality.b}
Let $M$ be a complete noncompact manifold of dimension $n+m$ with nonnegative sectional curvature. Let $\Sigma$ be a compact submanifold of $M$ of dimension $n$ (possibly with boundary $\partial \Sigma$), and let $f$ be a positive smooth function on $\Sigma$. If $m \geq 2$, then 
\[\int_\Sigma \sqrt{|\nabla^\Sigma f|^2 + f^2 \, |H|^2} + \int_{\partial \Sigma} f \geq n \, \Big ( \frac{(n+m) \, |B^{n+m}|}{m \, |B^m|} \Big )^{\frac{1}{n}} \, \theta^{\frac{1}{n}} \, \Big ( \int_\Sigma f^{\frac{n}{n-1}} \Big )^{\frac{n-1}{n}},\] 
where $\theta$ denotes the asymptotic volume ratio of $M$ and $H$ denotes the mean curvature vector of $\Sigma$.
\end{theorem}

Note that $(n+2) \, |B^{n+2}| = 2 \, |B^2| \, |B^n|$. Hence, we obtain a sharp Sobolev inequality for submanifolds of codimension $2$: 

\begin{corollary}
\label{codim.2}
Let $M$ be a complete noncompact manifold of dimension $n+2$ with nonnegative sectional curvature. Let $\Sigma$ be a compact submanifold of $M$ of dimension $n$ (possibly with boundary $\partial \Sigma$), and let $f$ be a positive smooth function on $\Sigma$. Then 
\[\int_\Sigma \sqrt{|\nabla^\Sigma f|^2 + f^2 \, |H|^2} + \int_{\partial \Sigma} f \geq n \, |B^n|^{\frac{1}{n}} \, \theta^{\frac{1}{n}} \, \Big ( \int_\Sigma f^{\frac{n}{n-1}} \Big )^{\frac{n-1}{n}},\] 
where $\theta$ denotes the asymptotic volume ratio of $M$ and $H$ denotes the mean curvature vector of $\Sigma$.
\end{corollary}

Moreover, we can characterize the case of equality in Corollary \ref{codim.2}:

\begin{theorem}
\label{rigidity.theorem.b}
Let $M$ be a complete noncompact manifold of dimension $n+2$ with nonnegative sectional curvature. Let $\Sigma$ be a compact submanifold of $M$ of dimension $n$ (possibly with boundary $\partial \Sigma$), and let $f$ be a positive smooth function on $\Sigma$. Suppose that 
\[\int_\Sigma \sqrt{|\nabla^\Sigma f|^2 + f^2 \, |H|^2} + \int_{\partial \Sigma} f = n \, |B^n|^{\frac{1}{n}} \, \theta^{\frac{1}{n}} \, \Big ( \int_\Sigma f^{\frac{n}{n-1}} \Big )^{\frac{n-1}{n}} > 0,\] 
where $\theta$ denotes the asymptotic volume ratio of $M$ and $H$ denotes the mean curvature vector of $\Sigma$. Then $f$ is constant, $M$ is isometric to Euclidean space, and $\Sigma$ is a flat round ball.
\end{theorem}

By putting $f=1$ in Corollary \ref{codim.2}, we obtain an isoperimetric inequality for minimal submanifolds of codimension $2$, generalizing the result in \cite{Brendle}.

\begin{corollary}
\label{isop.inequality.b}
Let $M$ be a complete noncompact manifold of dimension $n+2$ with nonnegative sectional curvature. Let $\Sigma$ be a compact minimal submanifold of $M$ of dimension $n$ with boundary $\partial \Sigma$. Then 
\[|\partial \Sigma| \geq n \, |B^n|^{\frac{1}{n}} \, \theta^{\frac{1}{n}} \, |\Sigma|^{\frac{n-1}{n}},\] 
where $\theta$ denotes the asymptotic volume ratio of the ambient manifold $M$.
\end{corollary}

%

Finally, the inequalities in Corollary \ref{codim.2} and Corollary \ref{isop.inequality.b} also hold in the codimension $1$ setting. Indeed, if $\Sigma$ is an $n$-dimensional submanifold of an $(n+1)$-dimensional manifold $M$, then we can view $\Sigma$ as a submanifold of the $(n+2)$-dimensional manifold $M \times \mathbb{R}$. Note that the product $M \times \mathbb{R}$ has the same asymptotic volume ratio as $M$ itself. 

The proof of Theorem \ref{sobolev.inequality.b} will be presented in Section \ref{proof.of.sobolev.inequality.b}. This argument extends our earlier proof in the Euclidean case (cf. \cite{Brendle}), and relies on the Alexandrov-Bakelman-Pucci technique. Moreover, the proof shares some common features with the work of E.~Heintze and H.~Karcher \cite{Heintze-Karcher} concerning the volume of a tubular neighborhood of a submanifold. Finally, the proof of Theorem \ref{rigidity.theorem.b} will be discussed in Section \ref{proof.of.rigidity.theorem.b}.

\section{Proof of Theorem \ref{sobolev.inequality.a}}

\label{proof.of.sobolev.inequality.a}

Throughout this section, we assume that $(M,g)$ is a complete noncompact manifold of dimension $n$ with nonnegative Ricci curvature. Moreover, we assume that $D$ is a compact domain in $M$, and $f$ is a positive smooth function on $D$. Let $R$ denote the Riemann curvature tensor of $(M,g)$. 

It suffices to prove the assertion in the special case when $D$ is connected. By scaling, we may assume that 
\[\int_D |\nabla f| + \int_{\partial D} f = n \int_D f^{\frac{n}{n-1}}.\] 
Since $D$ is connected, we can find a function $u: D \to \mathbb{R}$ with the property that 
\[\text{\rm div}(f \, \nabla u) = n \, f^{\frac{n}{n-1}} - |\nabla f|\] 
on $D$ and $\langle \nabla u,\eta \rangle = 1$ at each point on $\partial D$. Here, $\eta$ denotes the outward-pointing unit normal to $\partial D$. Standard elliptic regularity theory implies that the function $u$ is of class $C^{2,\gamma}$ for each $0 < \gamma < 1$ (cf. \cite{Gilbarg-Trudinger}, Theorem 6.30).

We define 
\[U := \{x \in D \setminus \partial D: |\nabla u(x)| < 1\}.\] 
For each $r > 0$, we denote by $A_r$ the set of all points $\bar{x} \in U$ with the property that 
\[r \, u(x) + \frac{1}{2} \, d \big ( x,\exp_{\bar{x}}(r \, \nabla u(\bar{x})) \big )^2 \geq r \, u(\bar{x}) + \frac{1}{2} \, r^2 \, |\nabla u(\bar{x})|^2\] 
for all $x \in D$. Moreover, for each $r > 0$, we define a map $\Phi_r: D \to M$ by 
\[\Phi_r(x) = \exp_x(r \, \nabla u(x))\] 
for all $x \in D$. Note that the map $\Phi_r$ is of class $C^{1,\gamma}$ for each $0 < \gamma < 1$. \\

\begin{lemma} 
\label{Laplacian.a}
Assume that $x \in U$. Then $\Delta u(x) \leq n \, f(x)^{\frac{1}{n-1}}$. 
\end{lemma}

\textbf{Proof.} 
Using the inequality $|\nabla u(x)| < 1$ and the Cauchy-Schwarz inequality, we obtain 
\[-\langle \nabla f(x),\nabla u(x) \rangle \leq |\nabla f(x)|.\] 
Moreover, $\text{\rm div}(f \, \nabla u) = n \, f^{\frac{n}{n-1}} - |\nabla f|$ by definition of $u$. This implies 
\[f(x) \, \Delta u(x) = n \, f(x)^{\frac{n}{n-1}} - |\nabla f(x)| - \langle \nabla f(x),\nabla u(x) \rangle \leq n \, f(x)^{\frac{n}{n-1}}.\] 
From this, the assertion follows. \\

\begin{lemma} 
\label{Phi.surjective.a}
The set 
\[\{p \in M: \text{\rm $d(x,p) < r$ for all $x \in D$}\}\] 
is contained in the set 
\[\{\Phi_r(x): x \in A_r\}.\]
\end{lemma}

\textbf{Proof.} 
Fix a point $p \in M$ with the property that $d(x,p) < r$ for all $x \in D$. Since $\langle \nabla u,\eta \rangle = 1$ at each point on $\partial D$, the function $x \mapsto r \, u(x) + \frac{1}{2} \, d(x,p)^2$ cannot attain its minimum on the boundary of $D$. Let us fix a point $\bar{x} \in D \setminus \partial D$ where the function $x \mapsto r \, u(x) + \frac{1}{2} \, d(x,p)^2$ attains its minimum. Moreover, let $\bar{\gamma}: [0,r] \to M$ be a minimizing geodesic such that $\bar{\gamma}(0) = \bar{x}$ and $\bar{\gamma}(r) = p$. Clearly, $r \, |\bar{\gamma}'(0)| = d(\bar{x},p)$. For every smooth path $\gamma: [0,r] \to M$ satisfying $\gamma(0) \in D$ and $\gamma(r) = p$, we obtain 
\begin{align*} 
r \, u(\gamma(0)) + \frac{1}{2} \, r \int_0^r |\gamma'(t)|^2 \, dt 
&\geq r \, u(\gamma(0)) + \frac{1}{2} \, d(\gamma(0),p)^2 \\ 
&\geq r \, u(\bar{x}) + \frac{1}{2} \, d(\bar{x},p)^2 \\ 
&= r \, u(\bar{\gamma}(0)) + \frac{1}{2} \, r^2 \, |\bar{\gamma}'(0)|^2 \\ 
&= r \, u(\bar{\gamma}(0)) + \frac{1}{2} \, r \int_0^r |\bar{\gamma}'(t)|^2 \, dt. 
\end{align*} 
In other words, the path $\bar{\gamma}$ minimizes the functional $u(\gamma(0)) + \frac{1}{2} \int_0^r |\gamma'(t)|^2 \, dt$ among all smooth paths $\gamma: [0,r] \to M$ satisfying $\gamma(0) \in D$ and $\gamma(r) = p$. Hence, the formula for the first variation of energy implies 
\[\nabla u(\bar{x}) = \bar{\gamma}'(0).\] 
From this, we deduce that 
\[\Phi_r(\bar{x}) = \exp_{\bar{x}}(r \, \nabla u(\bar{x})) = \exp_{\bar{\gamma}(0)}(r \, \bar{\gamma}'(0)) = \bar{\gamma}(r) = p.\] 
Moreover, 
\[r \, |\nabla u(\bar{x})| = r \, |\bar{\gamma}'(0)| = d(\bar{x},p).\] 
By assumption, $d(\bar{x},p) < r$. This implies $|\nabla u(\bar{x})| < 1$. Therefore, $\bar{x} \in U$. Finally, for each point $x \in D$, we have 
\begin{align*} 
r \, u(x) + \frac{1}{2} \, d \big ( x,\exp_{\bar{x}}(r \, \nabla u(\bar{x})) \big )^2 
&= r \, u(x) + \frac{1}{2} \, d(x,p)^2 \\ 
&\geq r \, u(\bar{x}) + \frac{1}{2} \, d(\bar{x},p)^2 \\ 
&= r \, u(\bar{x}) + \frac{1}{2} \, r^2 \, |\nabla u(\bar{x})|^2. 
\end{align*} 
Thus, $\bar{x} \in A_r$. This completes the proof of Lemma \ref{Phi.surjective.a}. \\

\begin{lemma} 
\label{second.variation.a}
Assume that $\bar{x} \in A_r$, and let $\bar{\gamma}(t) := \exp_{\bar{x}}(t \, \nabla u(\bar{x}))$ for all $t \in [0,r]$. If $Z$ is a smooth vector field along $\bar{\gamma}$ satisfying $Z(r) = 0$, then 
\[(D^2 u)(Z(0),Z(0)) + \int_0^r \big ( |D_t Z(t)|^2 - R(\bar{\gamma}'(t),Z(t),\bar{\gamma}'(t),Z(t)) \big ) \, dt \geq 0.\]
\end{lemma}

\textbf{Proof.} 
Let us consider an arbitrary smooth path $\gamma: [0,r] \to M$ satisfying $\gamma(0) \in D$ and $\gamma(r) = \bar{\gamma}(r)$. Since $\bar{x} \in A_r$, we obtain 
\begin{align*} 
r \, u(\gamma(0)) + \frac{1}{2} \, r \int_0^r |\gamma'(t)|^2 \, dt 
&\geq r \, u(\gamma(0)) + \frac{1}{2} \, d(\gamma(0),\gamma(r))^2 \\ 
&= r \, u(\gamma(0)) + \frac{1}{2} \, d \big ( \gamma(0),\exp_{\bar{x}}(r \, \nabla u(\bar{x})) \big )^2 \\ 
&\geq r \, u(\bar{x}) + \frac{1}{2} \, r^2 \, |\nabla u(\bar{x})|^2 \\ 
&= r \, u(\bar{\gamma}(0)) + \frac{1}{2} \, r \int_0^r |\bar{\gamma}'(t)|^2 \, dt. 
\end{align*} 
In other words, the path $\bar{\gamma}$ minimizes the functional $u(\gamma(0)) + \frac{1}{2} \int_0^r |\gamma'(t)|^2 \, dt$ among all smooth paths $\gamma: [0,r] \to M$ satisfying $\gamma(0) \in D$ and $\gamma(r) = \bar{\gamma}(r)$. Hence, the assertion follows from the formula for the second variation of energy. \\

%

\begin{lemma} 
\label{no.conjugate.points.a}
Assume that $\bar{x} \in A_r$, and let $\bar{\gamma}(t) := \exp_{\bar{x}}(t \, \nabla u(\bar{x}))$ for all $t \in [0,r]$. Moreover, let $\{e_1,\hdots,e_n\}$ be an orthonormal basis of $T_{\bar{x}} M$. Suppose that $W$ is a Jacobi field along $\bar{\gamma}$ satisfying $\langle D_t W(0),e_j \rangle = (D^2 u)(W(0),e_j)$ for each $1 \leq j \leq n$. If $W(\tau) = 0$ for some $\tau \in (0,r)$, then $W$ vanishes identically.
\end{lemma}

\textbf{Proof.} 
Suppose that $W(\tau) = 0$ for some $\tau \in (0,r)$. By assumption, 
\[\langle D_t W(0),W(0) \rangle = (D^2 u)(W(0),W(0)).\] 
Since $W$ is a Jacobi field, we obtain  
\begin{align*} 
&\int_0^\tau \big ( |D_t W(t)|^2 - R(\bar{\gamma}'(t),W(t),\bar{\gamma}'(t),W(t)) \big ) \, dt \\ 
&= \langle D_t W(\tau),W(\tau) \rangle - \langle D_t W(0),W(0) \rangle \\ 
&= -(D^2 u)(W(0),W(0)). 
\end{align*}
Let us define a vector field $\tilde{W}$ along $\bar{\gamma}$ by 
\[\tilde{W}(t) = \begin{cases} W(t) & \text{\rm for $t \in [0,\tau]$} \\ 0 & \text{\rm for $t \in [\tau,r]$.} \end{cases}\] 
Clearly, $\tilde{W}(r) = 0$. Moreover, 
\begin{align*} 
&\int_0^r \big ( |D_t \tilde{W}(t)|^2 - R(\bar{\gamma}'(t),\tilde{W}(t),\bar{\gamma}'(t),\tilde{W}(t)) \big ) \, dt \\ 
&= -(D^2 u)(\tilde{W}(0),\tilde{W}(0)). 
\end{align*} 
Using Lemma \ref{second.variation.a}, we conclude that 
\begin{align*} 
&\int_0^r \big ( |D_t Z(t)|^2 - R(\bar{\gamma}'(t),Z(t),\bar{\gamma}'(t),Z(t)) \big ) \, dt \\ 
&\geq \int_0^r \big ( |D_t \tilde{W}(t)|^2 - R(\bar{\gamma}'(t),\tilde{W}(t),\bar{\gamma}'(t),\tilde{W}(t)) \big ) \, dt
\end{align*}
for every smooth vector field $Z$ along $\bar{\gamma}$ satisfying $Z(0) = \tilde{W}(0)$ and $Z(r) = \tilde{W}(r)$. By approximation, this inequality holds for every vector field $Z$ which is piecewise $C^1$ and satisfies $Z(0) = \tilde{W}(0)$ and $Z(r) = \tilde{W}(r)$. In other words, the vector field $\tilde{W}$ minimizes the index form among all vector fields which are piecewise $C^1$ and have the same boundary values as $\tilde{W}$. Consequently, $\tilde{W}$ must be of class $C^1$. This implies $D_t W(\tau) = 0$. Since $W(\tau) = 0$ and $D_t W(\tau) = 0$, standard uniqueness results for ODE imply that $W$ vanishes identically. \\

\begin{proposition} 
\label{monotonicity.a}
Assume that $x \in A_r$. Then the function 
\[t \mapsto (1+t \, f(x)^{\frac{1}{n-1}})^{-n} \, |\det D\Phi_t(x)|\] 
is monotone decreasing for $t \in (0,r)$.
\end{proposition}

\textbf{Proof.} 
Fix an arbitrary point $\bar{x} \in A_r$. Let $\{e_1,\hdots,e_n\}$ be an orthonormal basis of $T_{\bar{x}} M$, and let $(x_1,\hdots,x_n)$ be a system of geodesic normal coordinates around $\bar{x}$ such that $\frac{\partial}{\partial x_i} = e_i$ at $\bar{x}$. Let $\bar{\gamma}(t) := \exp_{\bar{x}}(t \, \nabla u(\bar{x}))$ for all $t \in [0,r]$. For each $1 \leq i \leq n$, we denote by $E_i(t)$ the parallel transport of $e_i$ along $\bar{\gamma}$. Moreover, for each $1 \leq i \leq n$, we denote by $X_i(t)$ the unique Jacobi field along $\bar{\gamma}$ satisfying $X_i(0) = e_i$ and 
\[\langle D_t X_i(0),e_j \rangle = (D^2 u)(e_i,e_j)\] 
for all $1 \leq j \leq n$. It follows from Lemma \ref{no.conjugate.points.a} that $X_1(t),\hdots,X_n(t)$ are linearly independent for each $t \in (0,r)$. 

Let us define an $n \times n$-matrix $P(t)$ by 
\[P_{ij}(t) = \langle X_i(t),E_j(t) \rangle\] 
for $1 \leq i,j \leq n$. Moreover, we define an $n \times n$-matrix $S(t)$ by 
\[S_{ij}(t) = R(\bar{\gamma}'(t),E_i(t),\bar{\gamma}'(t),E_j(t))\] 
for $1 \leq i,j \leq n$. Clearly, $S(t)$ is symmetric. Moreover, since $M$ has nonnegative Ricci curvature, we know that $\text{\rm tr}(S(t)) \geq 0$. Since the vector fields $X_1(t),\hdots,X_n(t)$ are Jacobi fields, we obtain 
\[P''(t) = -P(t) S(t).\] 
Moreover, 
\[P_{ij}(0) = \delta_{ij}\] 
and 
\[P_{ij}'(0) = (D^2 u)(e_i,e_j).\] 
In particular, the matrix $P'(0) P(0)^T$ is symmetric. Moreover, the matrix 
\[\frac{d}{dt} (P'(t) P(t)^T) = P''(t) P(t)^T + P'(t) P'(t)^T = -P(t) S(t) P(t)^T + P'(t) P'(t)^T\] 
is symmetric for each $t$. Thus, we conclude that the matrix $P'(t) P(t)^T$ is symmetric for each $t$. 

Since $X_1(t),\hdots,X_n(t)$ are linearly independent for each $t \in (0,r)$, the matrix $P(t)$ is invertible for each $t \in (0,r)$. Since $P'(t) P(t)^T$ is symmetric for each $t \in (0,r)$, it follows that the matrix $Q(t) := P(t)^{-1} P'(t)$ is symmetric for each $t \in (0,r)$. The matrix $Q(t)$ satisfies the Riccati equation 
\[Q'(t) = P(t)^{-1} P''(t) - P(t)^{-1} P'(t) P(t)^{-1} P'(t) = -S(t) - Q(t)^2\] 
for all $t \in (0,r)$. Moreover, since $Q(t)$ is symmetric, we obtain $\text{\rm tr}(Q(t)^2) \geq \frac{1}{n} \, \text{\rm tr}(Q(t))^2$ for all $t \in (0,r)$. Since $\text{\rm tr}(S(t)) \geq 0$, it follows that 
\[\frac{d}{dt} \text{\rm tr}(Q(t)) \leq -\text{\rm tr}(Q(t)^2) \leq -\frac{1}{n} \, \text{\rm tr}(Q(t))^2\] 
for all $t \in (0,r)$. Clearly, 
\[\lim_{t \to 0} Q_{ij}(t) = (D^2 u)(e_i,e_j).\]
Using Lemma \ref{Laplacian.a}, we obtain 
\[\lim_{t \to 0} \text{\rm tr}(Q(t)) = \Delta u(\bar{x}) \leq n \, f(\bar{x})^{\frac{1}{n-1}}.\] 
Hence, a standard ODE comparison principle implies 
\[\text{\rm tr}(Q(t)) \leq \frac{n \, f(\bar{x})^{\frac{1}{n-1}}}{1+t \, f(\bar{x})^{\frac{1}{n-1}}}\] 
for all $t \in (0,r)$.

We next consider the determinant of $P(t)$. Clearly, $\det P(t) > 0$ if $t$ is sufficiently small. Since $P(t)$ is invertible for each $t \in (0,r)$, it follows that $\det P(t) > 0$ for all $t \in (0,r)$. Using the estimate for the trace of $Q(t) = P(t)^{-1} P'(t)$, we obtain 
\[\frac{d}{dt} \log \det P(t) = \text{\rm tr}(Q(t)) \leq \frac{n \, f(\bar{x})^{\frac{1}{n-1}}}{1+t \, f(\bar{x})^{\frac{1}{n-1}}}\] 
for all $t \in (0,r)$. Consequently, the function 
\[t \mapsto (1+t \, f(\bar{x})^{\frac{1}{n-1}})^{-n} \, \det P(t)\] 
is monotone decreasing for $t \in (0,r)$. 

Finally, we observe that 
\[\frac{\partial \Phi_t}{\partial x_i}(\bar{x}) = X_i(t)\] 
for $1 \leq i \leq n$. Consequently, $|\det D\Phi_t(\bar{x})| = \det P(t)$ for all $t \in (0,r)$. Putting these facts together, the assertion follows. \\

\begin{corollary}
\label{upper.bound.for.Jacobian.a}
The Jacobian determinant of $\Phi_r$ satisfies 
\[|\det D\Phi_r(x)| \leq (1+r \, f(x)^{\frac{1}{n-1}})^n\] 
for all $x \in A_r$. 
\end{corollary}

\textbf{Proof.}
Since $\lim_{t \to 0} |\det D\Phi_t(x)| = 1$, the assertion follows from Proposition \ref{monotonicity.a}. \\

After these preparations, we now complete the proof of Theorem \ref{sobolev.inequality.a}. Using Lemma \ref{Phi.surjective.a} and Corollary \ref{upper.bound.for.Jacobian.a}, we obtain 
\begin{align*} 
&|\{p \in M: \text{\rm $d(x,p) < r$ for all $x \in D$}\}| \\ 
&\leq \int_{A_r} |\det D\Phi_r(x)| \, d\text{\rm vol}(x) \\ 
&\leq \int_U (1+r \, f(x)^{\frac{1}{n-1}})^n \, d\text{\rm vol}(x)
\end{align*}
for all $r > 0$. Finally, we divide by $r^n$ and send $r \to \infty$. This gives 
\[|B^n| \, \theta \leq \int_U f^{\frac{n}{n-1}} \leq \int_D f^{\frac{n}{n-1}}.\] 
Thus, 
\[\int_D |\nabla f| + \int_{\partial D} f = n \int_D f^{\frac{n}{n-1}} \geq n \, |B^n|^{\frac{1}{n}} \, \theta^{\frac{1}{n}} \, \Big ( \int_D f^{\frac{n}{n-1}} \Big )^{\frac{n-1}{n}}.\]
This completes the proof of Theorem \ref{sobolev.inequality.a}.

\section{Proof of Theorem \ref{rigidity.theorem.a}}

\label{proof.of.rigidity.theorem.a}

Let $(M,g)$ be a complete noncompact manifold of dimension $n$ with nonnegative Ricci curvature. Let $D$ be a compact domain in $M$ with boundary $\partial D$, and let $f$ be a positive smooth function on $D$ satisfying 
\[\int_D |\nabla f| + \int_{\partial D} f = n \, |B^n|^{\frac{1}{n}} \, \theta^{\frac{1}{n}} \, \Big ( \int_D f^{\frac{n}{n-1}} \Big )^{\frac{n-1}{n}} > 0,\] 
where $\theta$ denotes the asymptotic volume ratio of $M$. 

If $D$ is disconnected, we may apply Theorem \ref{sobolev.inequality.a} to each connected component of $D$, and take the sum over all connected components. This will lead to a contradiction. Therefore, $D$ must be connected.

By scaling, we may assume that 
\[\int_D |\nabla f| + \int_{\partial D} f = n \, |B^n| \, \theta\] 
and 
\[\int_D f^{\frac{n}{n-1}} = |B^n| \, \theta.\] 
In particular, 
\[\int_D |\nabla f| + \int_{\partial D} f = n \int_D f^{\frac{n}{n-1}}.\] 
Since $D$ is connected, we can find a function $u: D \to \mathbb{R}$ such that 
\[\text{\rm div}(f \, \nabla u) = n \, f^{\frac{n}{n-1}} - |\nabla f|\] 
on $D$ and $\langle \nabla u,\eta \rangle = 1$ at each point on $\partial D$. Moreover, $u$ is of class $C^{2,\gamma}$ for each $0 < \gamma < 1$. Let us define $U$, $A_r$, and $\Phi_r$ as in Section \ref{proof.of.sobolev.inequality.a}.

\begin{lemma}
\label{lower.bound.for.det.P.a}
Assume that $x \in U$. Then $|\det D\Phi_t(x)| \geq (1+t \, f(x)^{\frac{1}{n-1}})^n$ for each $t > 0$. 
\end{lemma}

\textbf{Proof.} 
Let us fix a point $\bar{x} \in U$. Suppose that $|\det D\Phi_t(\bar{x})| < (1+t \, f(\bar{x})^{\frac{1}{n-1}})^n$ for some $t > 0$. Let us fix a real number $\varepsilon \in (0,1)$ such that 
\[|\det D\Phi_t(\bar{x})| < (1-\varepsilon) \, (1+t \, f(\bar{x})^{\frac{1}{n-1}})^n.\] 
By continuity, we can find an open neighborhood $V$ of the point $\bar{x}$ such that 
\[|\det D\Phi_t(x)| \leq (1-\varepsilon) \, (1+t \, f(x)^{\frac{1}{n-1}})^n\] 
for all $x \in V$. Using Proposition \ref{monotonicity.a}, we conclude that 
\[|\det D\Phi_r(x)| \leq (1-\varepsilon) \, (1+r \, f(x)^{\frac{1}{n-1}})^n\] 
for all $r > t$ and all $x \in A_r \cap V$. Using this fact together with Lemma \ref{Phi.surjective.a} and Corollary \ref{upper.bound.for.Jacobian.a}, we obtain 
\begin{align*} 
&|\{p \in M: \text{\rm $d(x,p) < r$ for all $x \in D$}\}| \\ 
&\leq \int_{A_r} |\det D\Phi_r(x)| \, d\text{\rm vol}(x) \\ 
&\leq \int_U (1-\varepsilon \cdot 1_V(x)) \, (1+r \, f(x)^{\frac{1}{n-1}})^n \, d\text{\rm vol}(x)
\end{align*} 
for all $r > t$. Finally, we divide by $r^n$ and send $r \to \infty$. This gives 
\[|B^n| \, \theta \leq \int_U (1-\varepsilon \cdot 1_V) \, f^{\frac{n}{n-1}} < \int_D f^{\frac{n}{n-1}} = |B^n| \, \theta.\] 
This is a contradiction. \\

\begin{lemma} 
\label{Hessian.of.u.a}
Assume that $x \in U$. Then $D^2 u(x) = f(x)^{\frac{1}{n-1}} \, g$. 
\end{lemma}

\textbf{Proof.}
Let us fix a point $\bar{x} \in U$. Let $\{e_1,\hdots,e_n\}$ be an orthonormal basis of $T_{\bar{x}} M$. We define $\bar{\gamma}(t) := \exp_{\bar{x}}(t \, \nabla u(\bar{x}))$ for all $t \geq 0$. For each $1 \leq i \leq n$, we denote by $E_i(t)$ the parallel transport of $e_i$ along $\bar{\gamma}$. Moreover, for each $1 \leq i \leq n$, we denote by $X_i(t)$ the unique Jacobi field along $\bar{\gamma}$ satisfying $X_i(0) = e_i$ and 
\[\langle D_t X_i(0),e_j \rangle = (D^2 u)(e_i,e_j)\] 
for all $1 \leq j \leq n$. Finally, we define an $n \times n$-matrix $P(t)$ by 
\[P_{ij}(t) = \langle X_i(t),E_j(t) \rangle\] 
for $1 \leq i,j \leq n$. 

By Lemma \ref{lower.bound.for.det.P.a}, we know that $|\det P(t)| \geq (1+t \, f(\bar{x})^{\frac{1}{n-1}})^n$ for all $t > 0$. Since $\det P(t) > 0$ if $t > 0$ is sufficiently small, we conclude that 
\[\det P(t) \geq (1+t \, f(\bar{x})^{\frac{1}{n-1}})^n\] 
for all $t > 0$. In particular, $P(t)$ is invertible for each $t > 0$.

We next define $Q(t) := P(t)^{-1} P'(t)$ for all $t > 0$. As in Section \ref{proof.of.sobolev.inequality.a}, we can show that the matrix $Q(t)$ is symmetric for each $t > 0$. The Riccati equation for $Q(t)$ gives 
\[\frac{d}{dt} \text{\rm tr}(Q(t)) \leq -\text{\rm tr}(Q(t)^2) \leq -\frac{1}{n} \, \text{\rm tr}(Q(t))^2\] 
for all $t > 0$. Moreover, 
\[\lim_{t \to 0} \text{\rm tr}(Q(t)) = \Delta u(\bar{x}) \leq n \, f(\bar{x})^{\frac{1}{n-1}}\] 
by Lemma \ref{Laplacian.a}. This implies 
\[\text{\rm tr}(Q(t)) \leq \frac{n \, f(\bar{x})^{\frac{1}{n-1}}}{1+t \, f(\bar{x})^{\frac{1}{n-1}}},\] 
hence 
\[\frac{d}{dt} \log \det P(t) \leq \frac{n \, f(\bar{x})^{\frac{1}{n-1}}}{1+t \, f(\bar{x})^{\frac{1}{n-1}}}\] 
for all $t > 0$. Integrating this ODE gives 
\[\det P(t) \leq (1+t \, f(\bar{x})^{\frac{1}{n-1}})^n\] 
for all $t > 0$.

Putting these facts together, we conclude that $\det P(t) = (1+t \, f(\bar{x})^{\frac{1}{n-1}})^n$ for all $t > 0$. Differentiating this identity with respect to $t$, we obtain 
\[\text{\rm tr}(Q(t)) = \frac{n \, f(\bar{x})^{\frac{1}{n-1}}}{1+t \, f(\bar{x})^{\frac{1}{n-1}}}\] 
for all $t > 0$. Using the Riccati equation for $Q(t)$, we conclude that $\text{\rm tr}(Q(t)^2) = \frac{1}{n} \, \text{\rm tr}(Q(t))^2$ for all $t > 0$. Consequently, the trace-free part of $Q(t)$ vanishes for each $t > 0$. Therefore, 
\[Q_{ij}(t) = \frac{f(\bar{x})^{\frac{1}{n-1}}}{1+t \, f(\bar{x})^{\frac{1}{n-1}}} \, \delta_{ij}\] 
for all $t > 0$. In particular,  
\[(D^2 u)(e_i,e_j) = \lim_{t \to 0} Q_{ij}(t) = f(\bar{x})^{\frac{1}{n-1}} \, \delta_{ij}.\] 
This completes the proof of Lemma \ref{Hessian.of.u.a}. \\

\begin{lemma} 
\label{gradient.of.f.a}
Assume that $x \in U$. Then $\nabla f(x) = 0$. 
\end{lemma} 

\textbf{Proof.} 
Let us consider an arbitrary point $x \in U$. Using the definition of $u$, we obtain 
\[f(x) \, \Delta u(x) = n \, f(x)^{\frac{n}{n-1}} - |\nabla f(x)| - \langle \nabla f(x),\nabla u(x) \rangle.\] 
On the other hand, Lemma \ref{Hessian.of.u.a} implies $\Delta u(x) = n \, f(x)^{\frac{1}{n-1}}$. Putting these facts together, we conclude that $\langle \nabla f(x),\nabla u(x) \rangle = -|\nabla f(x)|$. Since $|\nabla u(x)| < 1$, it follows that $\nabla f(x) = 0$. This completes the proof of Lemma \ref{gradient.of.f.a}. \\

\begin{lemma}
\label{density.a}
The set $U$ is dense in $D$. 
\end{lemma}

\textbf{Proof.} 
Suppose that $U$ is not dense in $D$. Arguing as in Section \ref{proof.of.sobolev.inequality.a}, we obtain 
\[|B^n| \, \theta \leq \int_U f^{\frac{n}{n-1}} < \int_D f^{\frac{n}{n-1}} = |B^n| \, \theta.\] 
This is a contradiction. This completes the proof of Lemma \ref{density.a}. \\

Since $U$ is a dense subset of $D$, we conclude that $\nabla f = 0$ and $D^2 u = f^{\frac{1}{n-1}} \, g$ at each point on $D$. Since $D$ is connected, it follows that $f$ is constant. This implies $|\partial D| = n \, |B^n|^{\frac{1}{n}} \, \theta^{\frac{1}{n}} \, |D|^{\frac{n-1}{n}}$. 

Note that $u$ is a smooth function on $D$. Each critical point of $u$ lies in the interior of $D$ and is nondegenerate with Morse index $0$. In particular, the function $u$ has at most finitely many critical points. 

We next consider the flow on $D$ generated by the vector field $-\nabla u$. Since the vector field $-\nabla u$ points inward along the boundary $\partial D$, the flow is defined for all nonnegative times. This gives a one-parameter family of smooth maps $\psi_s: D \to D$, where $s \geq 0$. Since $D$ is connected, standard arguments from Morse theory imply that the function $u$ has exactly one critical point, and $u$ attains its global minimum at that point. It follows that the diameter of $\psi_s(D)$ converges to $0$ as $s \to \infty$.

Since $D^2 u$ is a constant multiple of the metric, the isoperimetric ratio is unchanged under the flow $\psi_s$. This implies 
\[|\psi_s(\partial D)| = n \, |B^n|^{\frac{1}{n}} \, \theta^{\frac{1}{n}} \, |\psi_s(D)|^{\frac{n-1}{n}}\] 
for each $s \geq 0$. If $\theta<1$, this contradicts the Euclidean isoperimetric inequality when $s$ is sufficiently large. Thus, we conclude that $\theta=1$. Consequently, $M$ is isometric to Euclidean space.

Once we know that $M$ is isometric to Euclidean space, it follows that $D$ is a round ball. This completes the proof of Theorem \ref{rigidity.theorem.a}.

\section{Proof of Theorem \ref{sobolev.inequality.b}}

\label{proof.of.sobolev.inequality.b}

Throughout this section, we assume that $(M,\bar{g})$ is a complete noncompact manifold of dimension $n+m$ with nonnegative sectional curvature. Moreover, we assume that $\Sigma$ is a compact submanifold of $M$ of dimension $n$ (possibly with boundary $\partial \Sigma$), and $f$ is a positive smooth function on $\Sigma$. Let $\bar{D}$ denote the Levi-Civita connection on the ambient manifold $(M,\bar{g})$, and let $\bar{R}$ denote the Riemann curvature tensor of $(M,\bar{g})$. We denote by $I\!I$ the second fundamental form of $\Sigma$. For each point $x \in \Sigma$, $I\!I$ is a symmetric bilinear form on $T_x \Sigma$ which takes values in the normal space $T_x^\perp \Sigma$. If $X$ and $Y$ are tangent vector fields on $\Sigma$ and $V$ is a normal vector field along $\Sigma$, then $\langle I\!I(X,Y),V \rangle = \langle \bar{D}_X Y,V \rangle = -\langle \bar{D}_X V,Y \rangle$.

It suffices to prove the assertion in the special case when $\Sigma$ is connected. By scaling, we may assume that 
\[\int_\Sigma \sqrt{|\nabla^\Sigma f|^2 + f^2 \, |H|^2} + \int_{\partial \Sigma} f = n \int_\Sigma f^{\frac{n}{n-1}}.\] 
Since $\Sigma$ is connected, we can find a function $u: \Sigma \to \mathbb{R}$ with the property that 
\[\text{\rm div}_\Sigma(f \, \nabla^\Sigma u) = n \, f^{\frac{n}{n-1}} - \sqrt{|\nabla^\Sigma f|^2 + f^2 \, |H|^2}\] 
on $\Sigma$ and $\langle \nabla^\Sigma u,\eta \rangle = 1$ at each point on $\partial \Sigma$. Here, $\eta$ denotes the co-normal to $\partial \Sigma$. Standard elliptic regularity theory implies that the function $u$ is of class $C^{2,\gamma}$ for each $0 < \gamma < 1$ (cf. \cite{Gilbarg-Trudinger}, Theorem 6.30).

We define 
\begin{align*} 
\Omega &:= \{x \in \Sigma \setminus \partial \Sigma: |\nabla^\Sigma u(x)| < 1\}, \\ 
U &:= \{(x,y): x \in \Sigma \setminus \partial \Sigma, \, y \in T_x^\perp \Sigma, \, |\nabla^\Sigma u(x)|^2 + |y|^2 < 1\}. 
\end{align*} 
For each $r > 0$, we denote by $A_r$ the set of all points $(\bar{x},\bar{y}) \in U$ with the property that 
\[r \, u(x) + \frac{1}{2} \, d \big ( x,\exp_{\bar{x}}(r \, \nabla^\Sigma u(\bar{x}) + r \, \bar{y}) \big )^2 \geq r \, u(\bar{x}) + \frac{1}{2} \, r^2 \, (|\nabla^\Sigma u(\bar{x})|^2+|\bar{y}|^2)\] 
for all $x \in \Sigma$. Moreover, for each $r > 0$, we define a map $\Phi_r: T^\perp \Sigma \to M$ by 
\[\Phi_r(x,y) = \exp_x(r \, \nabla^\Sigma u(x) + r \, y)\] 
for all $x \in \Sigma$ and $y \in T_x^\perp \Sigma$. Note that the map $\Phi_r$ is of class $C^{1,\gamma}$ for each $0 < \gamma < 1$. \\

\begin{lemma} 
\label{Laplacian.b}
Assume that $x \in \Omega$ and $y \in T_x^\perp \Sigma$ satisfy $|\nabla^\Sigma u(x)|^2 + |y|^2 \leq 1$. Then $\Delta_\Sigma u(x) - \langle H(x),y \rangle \leq n \, f(x)^{\frac{1}{n-1}}$. 
\end{lemma}

\textbf{Proof.} 
Using the inequality $|\nabla^\Sigma u(x)|^2+|y|^2 \leq 1$ and the Cauchy-Schwarz inequality, we obtain 
\begin{align*} 
&-\langle \nabla^\Sigma f(x),\nabla^\Sigma u(x) \rangle - f(x) \, \langle H(x),y \rangle \\ 
&\leq \sqrt{|\nabla^\Sigma f(x)|^2 + f(x)^2 \, |H(x)|^2} \, \sqrt{|\nabla^\Sigma u(x)|^2+|y|^2} \\ 
&\leq \sqrt{|\nabla^\Sigma f(x)|^2 + f(x)^2 \, |H(x)|^2}. 
\end{align*}
Moreover, $\text{\rm div}_\Sigma(f \, \nabla^\Sigma u) = n \, f^{\frac{n}{n-1}} - \sqrt{|\nabla^\Sigma f|^2 + f^2 \, |H|^2}$ by definition of $u$. Consequently, 
\begin{align*} 
f(x) \, \Delta_\Sigma u(x) - f(x) \, \langle H(x),y \rangle 
&= n \, f(x)^{\frac{n}{n-1}} - \sqrt{|\nabla^\Sigma f(x)|^2 + f(x)^2 \, |H(x)|^2} \\ 
&- \langle \nabla^\Sigma f(x),\nabla^\Sigma u(x) \rangle - f(x) \, \langle H(x),y \rangle \\ 
&\leq n \, f(x)^{\frac{n}{n-1}}. 
\end{align*} 
From this, the assertion follows. \\

\begin{lemma} 
\label{Phi.surjective.b}
For each $0 \leq \sigma < 1$, the set 
\[\{p \in M: \text{\rm $\sigma r < d(x,p) < r$ for all $x \in \Sigma$}\}\] 
is contained in the set 
\[\{\Phi_r(x,y): (x,y) \in A_r, \, |\nabla^\Sigma u(x)|^2+|y|^2 > \sigma^2\}.\]
\end{lemma}

\textbf{Proof.} 
Let us fix a real number $0 \leq \sigma < 1$ and a point $p \in M$ with the property that $\sigma r < d(x,p) < r$ for all $x \in \Sigma$. Since $\langle \nabla^\Sigma u,\eta \rangle = 1$ at each point on $\partial \Sigma$, the function $x \mapsto r \, u(x) + \frac{1}{2} \, d(x,p)^2$ cannot attain its minimum on the boundary of $\Sigma$. Let us fix a point $\bar{x} \in \Sigma \setminus \partial \Sigma$ where the function $x \mapsto r \, u(x) + \frac{1}{2} \, d(x,p)^2$ attains its minimum. Moreover, let $\bar{\gamma}: [0,r] \to M$ be a minimizing geodesic such that $\bar{\gamma}(0) = \bar{x}$ and $\bar{\gamma}(r) = p$. Clearly, $r \, |\bar{\gamma}'(0)| = d(\bar{x},p)$. For every smooth path $\gamma: [0,r] \to M$ satisfying $\gamma(0) \in \Sigma$ and $\gamma(r) = p$, we obtain 
\begin{align*} 
r \, u(\gamma(0)) + \frac{1}{2} \, r \int_0^r |\gamma'(t)|^2 \, dt 
&\geq r \, u(\gamma(0)) + \frac{1}{2} \, d(\gamma(0),p)^2 \\ 
&\geq r \, u(\bar{x}) + \frac{1}{2} \, d(\bar{x},p)^2 \\ 
&= r \, u(\bar{\gamma}(0)) + \frac{1}{2} \, r^2 \, |\bar{\gamma}'(0)|^2 \\ 
&= r \, u(\bar{\gamma}(0)) + \frac{1}{2} \, r \int_0^r |\bar{\gamma}'(t)|^2 \, dt. 
\end{align*} 
In other words, the path $\bar{\gamma}$ minimizes the functional $u(\gamma(0)) + \frac{1}{2} \int_0^r |\gamma'(t)|^2 \, dt$ among all smooth paths $\gamma: [0,r] \to M$ satisfying $\gamma(0) \in \Sigma$ and $\gamma(r) = p$. Hence, the formula for the first variation of energy implies 
\[\nabla^\Sigma u(\bar{x}) - \bar{\gamma}'(0) \in T_{\bar{x}}^\perp \Sigma.\] 
Consequently, we can find a vector $\bar{y} \in T_{\bar{x}}^\perp \Sigma$ such that 
\[\nabla^\Sigma u(\bar{x}) + \bar{y} = \bar{\gamma}'(0).\] 
From this, we deduce that 
\[\Phi_r(\bar{x},\bar{y}) = \exp_{\bar{x}}(r \, \nabla^\Sigma u(\bar{x}) + r \, \bar{y}) = \exp_{\bar{\gamma}(0)}(r \, \bar{\gamma}'(0)) = \bar{\gamma}(r) = p.\] 
Moreover, 
\[r^2 \, (|\nabla^\Sigma u(\bar{x})|^2+|\bar{y}|^2) = r^2 \, |\nabla^\Sigma u(\bar{x})+\bar{y}|^2 = r^2 \, |\bar{\gamma}'(0)|^2 = d(\bar{x},p)^2.\] 
By assumption, $\sigma r < d(\bar{x},p) < r$. This implies $\sigma^2 < |\nabla^\Sigma u(\bar{x})|^2+|\bar{y}|^2 < 1$. In particular, $(\bar{x},\bar{y}) \in U$. Finally, for each point $x \in \Sigma$, we have 
\begin{align*} 
r \, u(x) + \frac{1}{2} \, d \big ( x,\exp_{\bar{x}}(r \, \nabla^\Sigma u(\bar{x}) + r \, \bar{y}) \big )^2 
&= r \, u(x) + \frac{1}{2} \, d(x,p)^2 \\ 
&\geq r \, u(\bar{x}) + \frac{1}{2} \, d(\bar{x},p)^2 \\ 
&= r \, u(\bar{x}) + \frac{1}{2} \, r^2 \, (|\nabla^\Sigma u(\bar{x})|^2+|\bar{y}|^2). 
\end{align*} 
Thus, $(\bar{x},\bar{y}) \in A_r$. This completes the proof of Lemma \ref{Phi.surjective.b}. \\

\begin{lemma} 
\label{second.variation.b}
Assume that $(\bar{x},\bar{y}) \in A_r$, and let $\bar{\gamma}(t) := \exp_{\bar{x}}(t \, \nabla^\Sigma u(\bar{x}) + t \, \bar{y})$ for all $t \in [0,r]$. If $Z$ is a smooth vector field along $\bar{\gamma}$ satisfying $Z(0) \in T_{\bar{x}} \Sigma$ and $Z(r) = 0$, then 
\begin{align*} 
&(D_\Sigma^2 u)(Z(0),Z(0)) - \langle I\!I(Z(0),Z(0)),\bar{y} \rangle \\ 
&+ \int_0^r \big ( |\bar{D}_t Z(t)|^2 - \bar{R}(\bar{\gamma}'(t),Z(t),\bar{\gamma}'(t),Z(t)) \big ) \, dt \geq 0. 
\end{align*}
\end{lemma}

\textbf{Proof.} 
Let us consider an arbitrary smooth path $\gamma: [0,r] \to M$ satisfying $\gamma(0) \in \Sigma$ and $\gamma(r) = \bar{\gamma}(r)$. Since $(\bar{x},\bar{y}) \in A_r$, we obtain 
\begin{align*} 
r \, u(\gamma(0)) + \frac{1}{2} \, r \int_0^r |\gamma'(t)|^2 \, dt 
&\geq r \, u(\gamma(0)) + \frac{1}{2} \, d(\gamma(0),\gamma(r))^2 \\ 
&= r \, u(\gamma(0)) + \frac{1}{2} \, d \big ( \gamma(0),\exp_{\bar{x}}(r \, \nabla^\Sigma u(\bar{x}) + r \, \bar{y}) \big )^2 \\ 
&\geq r \, u(\bar{x}) + \frac{1}{2} \, r^2 \, (|\nabla^\Sigma u(\bar{x})|^2+|\bar{y}|^2) \\ 
&= r \, u(\bar{\gamma}(0)) + \frac{1}{2} \, r \int_0^r |\bar{\gamma}'(t)|^2 \, dt. 
\end{align*} 
In other words, the path $\bar{\gamma}$ minimizes the functional $u(\gamma(0)) + \frac{1}{2} \int_0^r |\gamma'(t)|^2 \, dt$ among all smooth paths $\gamma: [0,r] \to M$ satisfying $\gamma(0) \in \Sigma$ and $\gamma(r) = \bar{\gamma}(r)$. Using the formula for the second variation of energy, we obtain 
\begin{align*} 
&(D_\Sigma^2 u)(Z(0),Z(0)) - \langle I\!I(Z(0),Z(0)),\bar{\gamma}'(0) \rangle \\ 
&+ \int_0^r \big ( |\bar{D}_t Z(t)|^2 - \bar{R}(\bar{\gamma}'(t),Z(t),\bar{\gamma}'(t),Z(t)) \big ) \, dt \geq 0. 
\end{align*}
On the other hand, the identity $\bar{\gamma}'(0) = \nabla^\Sigma u(\bar{x}) + \bar{y}$ implies 
\[\langle I\!I(Z(0),Z(0)),\bar{\gamma}'(0) \rangle = \langle I\!I(Z(0),Z(0)),\bar{y} \rangle.\] 
Putting these facts together, the assertion follows. \\

\begin{lemma}
\label{positivity.b}
Assume that $(\bar{x},\bar{y}) \in A_r$. Then $g+r \, D_\Sigma^2 u(\bar{x}) - r \, \langle I\!I(\bar{x}),\bar{y} \rangle \geq 0$.
\end{lemma}

\textbf{Proof.} 
As above, we define $\bar{\gamma}(t) := \exp_{\bar{x}}(t \, \nabla u(\bar{x}) + t \, \bar{y})$ for all $t \in [0,r]$. Let us fix an arbitrary vector $w \in T_{\bar{x}} \Sigma$, and let $W(t)$ denote the parallel transport of $w$ along $\bar{\gamma}$. Applying Lemma \ref{second.variation.b} to the vector field $Z(t) := (r-t) \, W(t)$ gives 
\begin{align*} 
&r \, g(w,w) + r^2 \, (D_\Sigma^2 u)(w,w) - r^2 \, \langle I\!I(w,w),\bar{y} \rangle \\ 
&- \int_0^r (r-t)^2 \, \bar{R}(\bar{\gamma}'(t),W(t),\bar{\gamma}'(t),W(t)) \, dt \geq 0. 
\end{align*} 
Since $M$ has nonnegative sectional curvature, it follows that 
\[r \, g(w,w) + r^2 \, (D_\Sigma^2 u)(w,w) - r^2 \, \langle I\!I(w,w),\bar{y} \rangle \geq 0,\] 
as claimed. \\

\begin{lemma} 
\label{no.conjugate.points.b}
Assume that $(\bar{x},\bar{y}) \in A_r$, and let $\bar{\gamma}(t) := \exp_{\bar{x}}(t \, \nabla^\Sigma u(\bar{x}) + t \, \bar{y})$ for all $t \in [0,r]$. Moreover, let $\{e_1,\hdots,e_n\}$ be an orthonormal basis of $T_{\bar{x}} \Sigma$. Suppose that $W$ is a Jacobi field along $\bar{\gamma}$ satisfying $W(0) \in T_{\bar{x}} \Sigma$ and $\langle \bar{D}_t W(0),e_j \rangle = (D_\Sigma^2 u)(W(0),e_j) - \langle I\!I(W(0),e_j),\bar{y} \rangle$ for each $1 \leq j \leq n$. If $W(\tau) = 0$ for some $\tau \in (0,r)$, then $W$ vanishes identically.
\end{lemma}

\textbf{Proof.} 
Suppose that $W(\tau) = 0$ for some $\tau \in (0,r)$. By assumption, 
\[\langle \bar{D}_t W(0),W(0) \rangle = (D_\Sigma^2 u)(W(0),W(0)) - \langle I\!I(W(0),W(0)),\bar{y} \rangle.\] 
Since $W$ is a Jacobi field, we obtain  
\begin{align*} 
&\int_0^\tau \big ( |\bar{D}_t W(t)|^2 - \bar{R}(\bar{\gamma}'(t),W(t),\bar{\gamma}'(t),W(t)) \big ) \, dt \\ 
&= \langle \bar{D}_t W(\tau),W(\tau) \rangle - \langle \bar{D}_t W(0),W(0) \rangle \\ 
&= -(D_\Sigma^2 u)(W(0),W(0)) + \langle I\!I(W(0),W(0)),\bar{y} \rangle. 
\end{align*}
Let us define a vector field $\tilde{W}$ along $\bar{\gamma}$ by 
\[\tilde{W}(t) = \begin{cases} W(t) & \text{\rm for $t \in [0,\tau]$} \\ 0 & \text{\rm for $t \in [\tau,r]$.} \end{cases}\] 
Clearly, $\tilde{W}(0) = W(0) \in T_{\bar{x}} \Sigma$ and $\tilde{W}(r) = 0$. Moreover, 
\begin{align*} 
&\int_0^r \big ( |\bar{D}_t \tilde{W}(t)|^2 - \bar{R}(\bar{\gamma}'(t),\tilde{W}(t),\bar{\gamma}'(t),\tilde{W}(t)) \big ) \, dt \\ 
&= -(D_\Sigma^2 u)(\tilde{W}(0),\tilde{W}(0)) + \langle I\!I(\tilde{W}(0),\tilde{W}(0)),\bar{y} \rangle. 
\end{align*}
Using Lemma \ref{second.variation.b}, we conclude that 
\begin{align*}
&\int_0^r \big ( |\bar{D}_t Z(t)|^2 - \bar{R}(\bar{\gamma}'(t),Z(t),\bar{\gamma}'(t),Z(t)) \big ) \, dt \\ 
&\geq \int_0^r \big ( |\bar{D}_t \tilde{W}(t)|^2 - \bar{R}(\bar{\gamma}'(t),\tilde{W}(t),\bar{\gamma}'(t),\tilde{W}(t)) \big ) \, dt 
\end{align*} 
for every smooth vector field $Z$ along $\bar{\gamma}$ satisfying $Z(0)=\tilde{W}(0)$ and $Z(r)=\tilde{W}(r)$. By approximation, this inequality holds for every vector field $Z$ which is piecewise $C^1$ and satisfies $Z(0) = \tilde{W}(0)$ and $Z(r) = \tilde{W}(r)$. In other words, the vector field $\tilde{W}$ minimizes the index form among all vector fields which are piecewise $C^1$ and have the same boundary values as $\tilde{W}$. Consequently, $\tilde{W}$ must be of class $C^1$. This implies $\bar{D}_t W(\tau) = 0$. Since $W(\tau) = 0$ and $\bar{D}_t W(\tau) = 0$, standard uniqueness results for ODE imply that $W$ vanishes identically. \\

\begin{proposition} 
\label{monotonicity.b}
Assume that $(x,y) \in A_r$. Then the function
\[t \mapsto t^{-m} \, (1+t \, f(x)^{\frac{1}{n-1}})^{-n} \, |\det D\Phi_t(x,y)|\] 
is monotone decreasing for $t \in (0,r)$.
\end{proposition}

\textbf{Proof.} 
Fix an arbitrary point $(\bar{x},\bar{y}) \in A_r$. Let us choose an orthonormal basis $\{e_1,\hdots,e_n\}$ of $T_{\bar{x}} M$ such that the $n \times n$-matrix 
\[(D_\Sigma^2 u)(e_i,e_j) - \langle I\!I(e_i,e_j),\bar{y} \rangle\] 
is diagonal. Let $(x_1,\hdots,x_n)$ be a system of geodesic normal coordinates on $\Sigma$ around the point $\bar{x}$. We can arrange that $\frac{\partial}{\partial x_i} = e_i$ at $\bar{x}$. Let $\{\nu_{n+1},\hdots,\nu_{n+m}\}$ be a local orthonormal frame for the normal bundle, chosen so that $\langle \bar{D}_{e_i} \nu_\alpha,\nu_\beta \rangle = 0$ at $\bar{x}$. We write a normal vector $y$ as $y = \sum_{\alpha=n+1}^{n+m} y_\alpha \nu_\alpha$. With this understood, $(x_1,\hdots,x_n,y_{n+1},\hdots,y_{n+m})$ is a local coordinate system on the total space of the normal bundle $T^\perp \Sigma$.

Let $\bar{\gamma}(t) := \exp_{\bar{x}}(t \, \nabla^\Sigma u(\bar{x}) + t \, \bar{y})$ for all $t \in [0,r]$. For each $1 \leq i \leq n$, we denote by $E_i(t)$ the parallel transport of $e_i$ along $\bar{\gamma}$. Moreover, for each $1 \leq i \leq n$, we denote by $X_i(t)$ the unique Jacobi field along $\bar{\gamma}$ satisfying $X_i(0) = e_i$ and 
\begin{align*} 
&\langle \bar{D}_t X_i(0),e_j \rangle = (D_\Sigma^2 u)(e_i,e_j) - \langle I\!I(e_i,e_j),\bar{y} \rangle, \\ 
&\langle \bar{D}_t X_i(0),\nu_\beta \rangle = \langle I\!I(e_i,\nabla^\Sigma u),\nu_\beta \rangle 
\end{align*} 
for all $1 \leq j \leq n$ and all $n+1 \leq \beta \leq n+m$. For each $n+1 \leq \alpha \leq n+m$, we denote by $N_\alpha(t)$ the parallel transport of $\nu_\alpha$ along $\bar{\gamma}$. Moreover, for each $n+1 \leq \alpha \leq n+m$, we denote by $Y_\alpha(t)$ the unique Jacobi field along $\bar{\gamma}$ satisfying $Y_\alpha(0) = 0$ and $\bar{D}_t Y_\alpha(0) = \nu_\alpha$. It follows from Lemma \ref{no.conjugate.points.b} that $X_1(t),\hdots,X_n(t),Y_{n+1},\hdots,Y_{n+m}(t)$ are linearly independent for each $t \in (0,r)$. 

Let us define an $(n+m) \times (n+m)$-matrix $P(t)$ by 
\[\begin{array}{ll}
P_{ij}(t) = \langle X_i(t),E_j(t) \rangle, & P_{i\beta}(t) = \langle X_i(t),N_\beta(t) \rangle, \\ 
P_{\alpha j}(t) = \langle Y_\alpha(t),E_i(t) \rangle, & P_{\alpha\beta}(t) = \langle Y_\alpha(t),N_\beta(t) \rangle 
\end{array}\] 
for $1 \leq i,j \leq n$ and $n+1 \leq \alpha,\beta \leq n+m$. Moreover, we define an $(n+m) \times (n+m)$-matrix $S(t)$ by 
\[\begin{array}{ll}
S_{ij}(t) = \bar{R}(\bar{\gamma}'(t),E_i(t),\bar{\gamma}'(t),E_j(t)), & S_{i\beta}(t) = \bar{R}(\bar{\gamma}'(t),E_i(t),\bar{\gamma}'(t),N_\beta(t)), \\ 
S_{\alpha j}(t) = \bar{R}(\bar{\gamma}'(t),N_\alpha(t),\bar{\gamma}'(t),E_j(t)), & S_{\alpha\beta}(t) = \bar{R}(\bar{\gamma}'(t),N_\alpha(t),\bar{\gamma}'(t),N_\beta(t)) 
\end{array}\] 
for $1 \leq i,j \leq n$ and $n+1 \leq \alpha,\beta \leq n+m$. Clearly, $S(t)$ is symmetric. Moreover, $S(t) \geq 0$ since $M$ has nonnegative sectional curvature. Since the vector fields $X_1(t),\hdots,X_n(t),Y_{n+1}(t),\hdots,Y_{n+m}(t)$ are Jacobi fields, we obtain 
\[P''(t) = -P(t) S(t).\] 
Moreover, 
\[P(0) = \begin{bmatrix} \delta_{ij}  & 0 \\ 0 & 0 \end{bmatrix}\] 
and 
\[P'(0) = \begin{bmatrix} (D_\Sigma^2 u)(e_i,e_j) - \langle I\!I(e_i,e_j),\bar{y} \rangle & \langle I\!I(e_i,\nabla^\Sigma u),\nu_\beta \rangle \\ 0 & \delta_{\alpha\beta} \end{bmatrix}.\] 
In particular, the matrix $P'(0) P(0)^T$ is symmetric. Moreover, the matrix 
\[\frac{d}{dt} (P'(t) P(t)^T) = P''(t) P(t)^T + P'(t) P'(t)^T = -P(t) S(t) P(t)^T + P'(t) P'(t)^T\] 
is symmetric for each $t$. Thus, we conclude that the matrix $P'(t) P(t)^T$ is symmetric for each $t$. 

Since $X_1(t),\hdots,X_n(t),Y_{n+1},\hdots,Y_{n+m}(t)$ are linearly independent for each $t \in (0,r)$, the matrix $P(t)$ is invertible for each $t \in (0,r)$. Since $P'(t) P(t)^T$ is symmetric for each $t \in (0,r)$, it follows that the matrix $Q(t) := P(t)^{-1} P'(t)$ is symmetric for each $t \in (0,r)$. The matrix $Q(t)$ satisfies the Riccati equation 
\[Q'(t) = P(t)^{-1} P''(t) - P(t)^{-1} P'(t) P(t)^{-1} P'(t) = -S(t) - Q(t)^2\] 
for all $t \in (0,r)$. Since $S(t) \geq 0$, it follows that 
\[Q'(t) \leq -Q(t)^2\] 
for all $t \in (0,r)$. Using the asymptotic expansion 
\[P(t) = \begin{bmatrix} \delta_{ij} + O(t) & O(t) \\ O(t) & t \, \delta_{\alpha\beta} + O(t^2) \end{bmatrix},\] 
we obtain 
\[P(t)^{-1} = \begin{bmatrix} \delta_{ij} + O(t) & O(1) \\ O(1) & t^{-1} \, \delta_{\alpha\beta} + O(1) \end{bmatrix}\] 
as $t \to 0$. 
Moreover, 
\[P'(t) = \begin{bmatrix} (D_\Sigma^2 u)(e_i,e_j) - \langle I\!I(e_i,e_j),\bar{y} \rangle + O(t) & O(1) \\ O(t) & \delta_{\alpha\beta} + O(t) \end{bmatrix}\] 
as $t \to 0$. Consequently, the matrix $Q(t) = P(t)^{-1} P'(t)$ satisfies the asymptotic expansion 
\[Q(t) = \begin{bmatrix} (D_\Sigma^2 u)(e_i,e_j) - \langle I\!I(e_i,e_j),\bar{y} \rangle + O(t) & O(1) \\ O(1) & t^{-1} \, \delta_{\alpha\beta} + O(1) \end{bmatrix}\]
as $t \to 0$. 

By our choice of $\{e_1,\hdots,e_n\}$, the matrix $(D_\Sigma^2 u)(e_i,e_j) - \langle I\!I(e_i,e_j),\bar{y} \rangle$ is diagonal. Let us write 
\[(D_\Sigma^2 u)(e_i,e_j) - \langle I\!I(e_i,e_j),\bar{y} \rangle = \lambda_i \, \delta_{ij}\] 
for $1 \leq i,j \leq n$. It follows from Lemma \ref{positivity.b} that $1+r\lambda_i \geq 0$ for each $1 \leq i \leq n$. Since 
\[Q(\tau) = \begin{bmatrix} \lambda_i \, \delta_{ij} + O(\tau) & O(1) \\ O(1) & \tau^{-1} \, \delta_{\alpha\beta} + O(1) \end{bmatrix}\] 
as $\tau \to 0$, we can find a small number $\tau_0 \in (0,r)$ such that 
\[Q(\tau) < \begin{bmatrix} (\lambda_i+\sqrt{\tau}) \, \delta_{ij} & 0 \\ 0 & 2\tau^{-1} \, \delta_{\alpha\beta} \end{bmatrix}\] 
for all $\tau \in (0,\tau_0)$. A standard ODE comparison principle implies 
\[Q(t) \leq \begin{bmatrix} \frac{(\lambda_i+\sqrt{\tau})}{1+(t-\tau)(\lambda_i+\sqrt{\tau})} \, \delta_{ij} & 0 \\ 0 & (t-\frac{\tau}{2})^{-1} \, \delta_{\alpha\beta} \end{bmatrix}\] 
for all $\tau \in (0,\tau_0)$ and all $t \in (\tau,r)$. Passing to the limit as $\tau \to 0$, we conclude that 
\[Q(t) \leq \begin{bmatrix} \frac{\lambda_i}{1+t\lambda_i} \, \delta_{ij} & 0 \\ 0 & t^{-1} \, \delta_{\alpha\beta} \end{bmatrix}\] 
for all $t \in (0,r)$. In particular, the trace of $Q(t)$ satisfies 
\[\text{\rm tr}(Q(t)) \leq \frac{m}{t} + \sum_{i=1}^n \frac{\lambda_i}{1+t\lambda_i}\] 
for all $t \in (0,r)$. It follows from Lemma \ref{Laplacian.b} that 
\[\sum_{i=1}^n \lambda_i = \Delta_\Sigma u(\bar{x}) - \langle H(\bar{x}),\bar{y} \rangle \leq n \, f(\bar{x})^{\frac{1}{n-1}}.\] 
Using the arithmetic-harmonic mean inequality, we obtain 
\[\sum_{i=1}^n \frac{1}{1+t\lambda_i} \geq \frac{n^2}{\sum_{i=1}^n (1+t\lambda_i)} \geq \frac{n}{1+t \, f(\bar{x})^{\frac{1}{n-1}}},\] 
hence 
\[\sum_{i=1}^n \frac{\lambda_i}{1+t\lambda_i} = \frac{1}{t} \, \Big ( n - \sum_{i=1}^n \frac{1}{1+t\lambda_i} \Big ) \leq \frac{n \, f(\bar{x})^{\frac{1}{n-1}}}{1+t \, f(\bar{x})^{\frac{1}{n-1}}}\] 
for all $t \in (0,r)$. Putting these facts together, we conclude that 
\[\text{\rm tr}(Q(t)) \leq \frac{m}{t} + \frac{n \, f(\bar{x})^{\frac{1}{n-1}}}{1+t \, f(\bar{x})^{\frac{1}{n-1}}}\] 
for all $t \in (0,r)$.

We next consider the determinant of $P(t)$. Clearly, $\lim_{t \to 0} t^{-m} \, \det P(t) = 1$. In particular, $\det P(t) > 0$ if $t>0$ is sufficiently small. Since $P(t)$ is invertible for each $t \in (0,r)$, it follows that $\det P(t) > 0$ for all $t \in (0,r)$. Using the estimate for the trace of $Q(t) = P(t)^{-1} P'(t)$, we obtain 
\[\frac{d}{dt} \log \det P(t) = \text{\rm tr}(Q(t)) \leq \frac{m}{t} + \frac{n \, f(\bar{x})^{\frac{1}{n-1}}}{1+t \, f(\bar{x})^{\frac{1}{n-1}}}\] 
for all $t \in (0,r)$. Consequently, the function 
\[t \mapsto t^{-m} \, (1+t \, f(\bar{x})^{\frac{1}{n-1}})^{-n} \, \det P(t)\] 
is monotone decreasing for $t \in (0,r)$. 

Finally, we observe that 
\[\frac{\partial \Phi_t}{\partial x_i}(\bar{x},\bar{y}) = X_i(t), \qquad \frac{\partial \Phi_t}{\partial y_\alpha}(\bar{x},\bar{y}) = Y_\alpha(t)\] 
for $1 \leq i \leq n$ and $n+1 \leq \alpha \leq n+m$. Consequently, $|\det D\Phi_t(\bar{x},\bar{y})| = \det P(t)$ for all $t \in (0,r)$. Putting these facts together, the assertion follows. \\

\begin{corollary} 
\label{upper.bound.for.Jacobian.b}
The Jacobian determinant of $\Phi_r$ satisfies 
\[|\det D\Phi_r(x,y)| \leq r^m \, (1+r \, f(x)^{\frac{1}{n-1}})^n\] 
for all $(x,y) \in A_r$. 
\end{corollary}

\textbf{Proof.} 
Since $\lim_{t \to 0} t^{-m} \, |\det D\Phi_t(x,y)| = 1$, the assertion follows from Proposition \ref{monotonicity.b}. \\

After these preparations, we now complete the proof of Theorem \ref{sobolev.inequality.b}. Using Lemma \ref{Phi.surjective.b} and Corollary \ref{upper.bound.for.Jacobian.b}, we obtain 
\begin{align*} 
&|\{p \in M: \text{\rm $\sigma r < d(x,p) < r$ for all $x \in \Sigma$}\}| \\ 
&\leq \int_\Omega \bigg ( \int_{\{y \in T_x^\perp \Sigma: \sigma^2 < |\nabla^\Sigma u(x)|^2+|y|^2 < 1\}} |\det D\Phi_r(x,y)| \, 1_{A_r}(x,y) \, dy \bigg ) \, d\text{\rm vol}(x) \\ 
&\leq \int_\Omega \bigg ( \int_{\{y \in T_x^\perp \Sigma: \sigma^2 < |\nabla^\Sigma u(x)|^2+|y|^2 < 1\}} r^m \, (1+r \, f(x)^{\frac{1}{n-1}})^n \, dy \bigg ) \, d\text{\rm vol}(x) \\ 
&= |B^m| \int_\Omega \Big [ (1-|\nabla^\Sigma u(x)|^2)^{\frac{m}{2}} - (\sigma^2 - |\nabla^\Sigma u(x)|^2)_+^{\frac{m}{2}} \Big ] \\ 
&\hspace{40mm} \cdot r^m \, (1+r \, f(x)^{\frac{1}{n-1}})^n \, d\text{\rm vol}(x)  
\end{align*}
for all $r > 0$ and all $0 \leq \sigma < 1$. Since $m \geq 2$, the mean value theorem implies $b^{\frac{m}{2}} - a^{\frac{m}{2}} \leq \frac{m}{2} \, (b-a)$ for $0 \leq a \leq b \leq 1$. Hence, we have the pointwise inequality 
\begin{align*} 
&(1-|\nabla^\Sigma u(x)|^2)^{\frac{m}{2}} - (\sigma^2 - |\nabla^\Sigma u(x)|^2)_+^{\frac{m}{2}} \\ 
&\leq \frac{m}{2} \, \Big [ (1-|\nabla^\Sigma u(x)|^2) -  (\sigma^2-|\nabla^\Sigma u(x)|^2)_+ \Big ] \leq \frac{m}{2} \, (1-\sigma^2) 
\end{align*}
for all $x \in \Omega$ and all $0 \leq \sigma < 1$. Therefore, 
\begin{align*} 
&|\{p \in M: \text{\rm $\sigma r < d(x,p) < r$ for all $x \in \Sigma$}\}| \\ 
&\leq \frac{m}{2} \, |B^m| \, (1-\sigma^2) \int_\Omega r^m \, (1+ r \, f(x)^{\frac{1}{n-1}})^n \, d\text{\rm vol}(x)
\end{align*} 
for all $r > 0$ and all $0 \leq \sigma < 1$. In the next step, we divide by $r^{n+m}$ and send $r \to \infty$ while keeping $\sigma$ fixed. This gives 
\[|B^{n+m}| \, (1-\sigma^{n+m}) \, \theta \leq \frac{m}{2} \, |B^m| \, (1-\sigma^2) \int_\Omega f^{\frac{n}{n-1}}\] 
for all $0 \leq \sigma < 1$. Finally, if we divide by $1-\sigma$ and send $\sigma \to 1$, we obtain 
\[(n+m) \, |B^{n+m}| \, \theta \leq m \, |B^m| \int_\Omega f^{\frac{n}{n-1}} \leq m \, |B^m| \int_\Sigma f^{\frac{n}{n-1}}.\] 
Thus, 
\begin{align*} 
&\int_\Sigma \sqrt{|\nabla^\Sigma f|^2 + f^2 \, |H|^2} + \int_{\partial \Sigma} f \\ 
&= n \int_\Sigma f^{\frac{n}{n-1}} \geq n \, \Big ( \frac{(n+m) \, |B^{n+m}|}{m \, |B^m|} \Big )^{\frac{1}{n}} \, \theta^{\frac{1}{n}} \, \Big ( \int_\Sigma f^{\frac{n}{n-1}} \Big )^{\frac{n-1}{n}}. 
\end{align*} 
This completes the proof of Theorem \ref{sobolev.inequality.b}.

\section{Proof of Theorem \ref{rigidity.theorem.b}}

\label{proof.of.rigidity.theorem.b}

Let $(M,\bar{g})$ be a complete noncompact manifold of dimension $n+2$ with nonnegative sectional curvature. Let $\Sigma$ be a compact submanifold of $M$ of dimension $n$ (possibly with boundary $\partial \Sigma$), and let $f$ be a positive smooth function on $\Sigma$ satisfying 
\[\int_\Sigma \sqrt{|\nabla^\Sigma f|^2 + f^2 \, |H|^2} + \int_{\partial \Sigma} f = n \, |B^n|^{\frac{1}{n}} \, \theta^{\frac{1}{n}} \, \Big ( \int_\Sigma f^{\frac{n}{n-1}} \Big )^{\frac{n-1}{n}} > 0,\] 
where $\theta$ denotes the asymptotic volume ratio of $M$. 

If $\Sigma$ is disconnected, we may apply Corollary \ref{codim.2} to each connected component of $\Sigma$, and take the sum over all connected components. This will lead to a contradiction. Therefore, $\Sigma$ must be connected.

By scaling, we may assume that 
\[\int_\Sigma \sqrt{|\nabla^\Sigma f|^2 + f^2 \, |H|^2} + \int_{\partial \Sigma} f = n \, |B^n| \, \theta\] 
and 
\[\int_\Sigma f^{\frac{n}{n-1}} = |B^n| \, \theta.\] 
In particular, 
\[\int_\Sigma \sqrt{|\nabla^\Sigma f|^2 + f^2 \, |H|^2} + \int_{\partial \Sigma} f = n \int_\Sigma f^{\frac{n}{n-1}}.\] 
Since $\Sigma$ is connected, we can find a function $u: \Sigma \to \mathbb{R}$ such that 
\[\text{\rm div}_\Sigma(f \, \nabla^\Sigma u) = n \, f^{\frac{n}{n-1}} - \sqrt{|\nabla^\Sigma f|^2 + f^2 \, |H|^2}\] 
on $\Sigma$ and $\langle \nabla^\Sigma u,\eta \rangle = 1$ at each point on $\partial \Sigma$. Moreover, $u$ is of class $C^{2,\gamma}$ for each $0 < \gamma < 1$. Let us define $\Omega$, $U$, $A_r$, and $\Phi_r$ as in Section \ref{proof.of.sobolev.inequality.b}.

\begin{lemma}
\label{lower.bound.for.det.P.b}
Assume that $x \in \Omega$ and $y \in T_x^\perp \Sigma$ satisfy $|\nabla^\Sigma u(x)|^2 + |y|^2 = 1$. Then $|\det D\Phi_t(x,y)| \geq t^2 \, (1+t \, f(x)^{\frac{1}{n-1}})^n$ for each $t > 0$. 
\end{lemma}

\textbf{Proof.} 
Let us fix a point $\bar{x} \in \Omega$ and a vector $\bar{y} \in T_{\bar{x}}^\perp \Sigma$ satisfying $|\nabla^\Sigma u(\bar{x})|^2 + |\bar{y}|^2 = 1$. Suppose that $|\det D\Phi_t(\bar{x},\bar{y})| < t^2 \, (1+t \, f(\bar{x})^{\frac{1}{n-1}})^n$ for some $t > 0$. Let us fix a real number $\varepsilon \in (0,1)$ such that 
\[|\det D\Phi_t(\bar{x},\bar{y})| < (1-\varepsilon) \, t^2 \, (1+t \, f(\bar{x})^{\frac{1}{n-1}})^n.\] 
By continuity, we can find an open neighborhood $V$ of the point $(\bar{x},\bar{y})$ such that 
\[|\det D\Phi_t(x,y)| \leq (1-\varepsilon) \, t^2 \, (1+t \, f(x)^{\frac{1}{n-1}})^n\] 
for all $(x,y) \in V$. Using Proposition \ref{monotonicity.b}, we conclude that 
\[|\det D\Phi_r(x,y)| \leq (1-\varepsilon) \, r^2 \, (1+r \, f(x)^{\frac{1}{n-1}})^n\] 
for all $r > t$ and all $(x,y) \in A_r \cap V$. Using this fact together with Lemma \ref{Phi.surjective.b} and Corollary \ref{upper.bound.for.Jacobian.b}, we obtain 
\begin{align*} 
&|\{p \in M: \text{\rm $\sigma r < d(x,p) < r$ for all $x \in \Sigma$}\}| \\ 
&\leq \int_\Omega \bigg ( \int_{\{y \in T_x^\perp \Sigma: \sigma^2 < |\nabla^\Sigma u(x)|^2+|y|^2 < 1\}} |\det D\Phi_r(x,y)| \, 1_{A_r}(x,y) \, dy \bigg ) \, d\text{\rm vol}(x) \\ 
&\leq \int_\Omega \bigg ( \int_{\{y \in T_x^\perp \Sigma: \sigma^2 < |\nabla^\Sigma u(x)|^2+|y|^2 < 1\}} (1-\varepsilon \cdot 1_V(x,y)) \\ 
&\hspace{60mm} \cdot r^2 \, (1+r \, f(x)^{\frac{1}{n-1}})^n \, dy \bigg ) \, d\text{\rm vol}(x) \\ 
&\leq |B^2| \, (1-\sigma^2) \int_\Omega r^2 \, (1+r \, f(x)^{\frac{1}{n-1}})^n \, d\text{\rm vol}(x) \\ 
&- \varepsilon \int_\Omega \bigg ( \int_{\{y \in T_x^\perp \Sigma: \sigma^2 < |\nabla^\Sigma u(x)|^2+|y|^2 < 1\}} 1_V(x,y) \, r^2 \, (1+r \, f(x)^{\frac{1}{n-1}})^n \, dy \bigg ) \, d\text{\rm vol}(x) 
\end{align*}
for all $r > t$ and all $0 \leq \sigma < 1$. We now divide by $r^{n+2}$ and send $r \to \infty$, while keeping $\sigma$ fixed. This implies 
\begin{align*} 
&|B^{n+2}| \, (1-\sigma^{n+2}) \, \theta \\ 
&\leq |B^2| \, (1-\sigma^2) \int_\Omega f(x)^{\frac{n}{n-1}} \, d\text{\rm vol}(x) \\ 
&- \varepsilon \int_\Omega \bigg ( \int_{\{y \in T_x^\perp \Sigma: \sigma^2 < |\nabla^\Sigma u(x)|^2+|y|^2 < 1\}} 1_V(x,y) \, f(x)^{\frac{n}{n-1}} \, dy \bigg ) \, d\text{\rm vol}(x) 
\end{align*} 
for all $0 \leq \sigma < 1$. Dividing by $1-\sigma$ and taking the limit as $\sigma \to 1$ gives 
\[(n+2) \, |B^{n+2}| \, \theta < 2 \, |B^2| \int_\Omega f^{\frac{n}{n-1}} \leq 2 \, |B^2| \, |B^n| \, \theta.\] 
This contradicts the fact that $(n+2) \, |B^{n+2}| = 2 \, |B^2| \, |B^n|$. \\

\begin{lemma} 
\label{modified.Hessian.of.u.b}
Assume that $x \in \Omega$ and $y \in T_x^\perp \Sigma$ satisfy $|\nabla^\Sigma u(x)|^2 + |y|^2 = 1$. Then $D_\Sigma^2 u(x) - \langle I\!I(x),y \rangle = f(x)^{\frac{1}{n-1}} \, g$. 
\end{lemma}

\textbf{Proof.} 
Let us fix a point $\bar{x} \in \Omega$ and a vector $\bar{y} \in T_{\bar{x}}^\perp \Sigma$ satisfying $|\nabla^\Sigma u(\bar{x})|^2+|\bar{y}|^2 = 1$. Let $\{e_1,\hdots,e_n\}$ be an orthonormal basis of $T_{\bar{x}} M$ with the property that the $n \times n$-matrix 
\[(D_\Sigma^2 u)(e_i,e_j) - \langle I\!I(e_i,e_j),\bar{y} \rangle\] 
is diagonal, and let $\{\nu_{n+1},\nu_{n+2}\}$ be an orthonormal basis of $T_{\bar{x}}^\perp \Sigma$. We define $\bar{\gamma}(t) := \exp_{\bar{x}}(t \, \nabla u(\bar{x}) + t \, \bar{y})$ for all $t \geq 0$. For each $1 \leq i \leq n$, we denote by $E_i(t)$ the parallel transport of $e_i$ along $\bar{\gamma}$. Moreover, for each $1 \leq i \leq n$, we denote by $X_i(t)$ the unique Jacobi field along $\bar{\gamma}$ satisfying $X_i(0) = e_i$ and 
\begin{align*} 
&\langle \bar{D}_t X_i(0),e_j \rangle = (D_\Sigma^2 u)(e_i,e_j) - \langle I\!I(e_i,e_j),\bar{y} \rangle, \\ 
&\langle \bar{D}_t X_i(0),\nu_\beta \rangle = \langle I\!I(e_i,\nabla^\Sigma u),\nu_\beta \rangle 
\end{align*} 
for all $1 \leq j \leq n$ and all $n+1 \leq \beta \leq n+2$. For each $n+1 \leq \alpha \leq n+2$, we denote by $N_\alpha(t)$ the parallel transport of $\nu_\alpha$ along $\bar{\gamma}$. Moreover, for each $n+1 \leq \alpha \leq n+2$, we denote by $Y_\alpha(t)$ the unique Jacobi field along $\bar{\gamma}$ satisfying $Y_\alpha(0) = 0$ and $\bar{D}_t Y_\alpha(0) = \nu_\alpha$. 

Finally, we define an $(n+2) \times (n+2)$-matrix $P(t)$ by 
\[\begin{array}{ll}
P_{ij}(t) = \langle X_i(t),E_j(t) \rangle, & P_{i\beta}(t) = \langle X_i(t),N_\beta(t) \rangle, \\ 
P_{\alpha j}(t) = \langle Y_\alpha(t),E_i(t) \rangle, & P_{\alpha\beta}(t) = \langle Y_\alpha(t),N_\beta(t) \rangle 
\end{array}\] 
for $1 \leq i,j \leq n$ and $n+1 \leq \alpha,\beta \leq n+2$. 

By Lemma \ref{lower.bound.for.det.P.b}, we know that $|\det P(t)| \geq t^2 \, (1+t \, f(\bar{x})^{\frac{1}{n-1}})^n$ for all $t > 0$. Since $\det P(t) > 0$ if $t > 0$ is sufficiently small, we conclude that 
\[\det P(t) \geq t^2 \, (1+t \, f(\bar{x})^{\frac{1}{n-1}})^n\] 
for all $t > 0$. In particular, $P(t)$ is invertible for each $t > 0$.

We next define $Q(t) := P(t)^{-1} P'(t)$ for all $t > 0$. Moreover, we write 
\[(D_\Sigma^2 u)(e_i,e_j) - \langle I\!I(e_i,e_j),\bar{y} \rangle = \lambda_i \, \delta_{ij}\] 
for $1 \leq i,j \leq n$. Arguing as in Section \ref{proof.of.sobolev.inequality.b}, we obtain 
\[\text{\rm tr}(Q(t)) \leq \frac{2}{t} + \sum_{i=1}^n \frac{\lambda_i}{1+t\lambda_i}\] 
for all $t>0$ satisfying $\min_{1 \leq i \leq n} (1+t\lambda_i) > 0$. Moreover, $\sum_{i=1}^n \lambda_i \leq n \, f(\bar{x})^{\frac{1}{n-1}}$ by Lemma \ref{Laplacian.b}. As above, the arithmetic-harmonic mean inequality implies 
\[\sum_{i=1}^n \frac{\lambda_i}{1+t\lambda_i} \leq \frac{n \, f(\bar{x})^{\frac{1}{n-1}}}{1+t \, f(\bar{x})^{\frac{1}{n-1}}}\] 
for all $t>0$ satisfying $\min_{1 \leq i \leq n} (1+t\lambda_i) > 0$. Therefore, 
\[\text{\rm tr}(Q(t)) \leq \frac{2}{t} + \frac{n \, f(\bar{x})^{\frac{1}{n-1}}}{1+t \, f(\bar{x})^{\frac{1}{n-1}}},\] 
hence 
\[\frac{d}{dt} \log \det P(t) \leq \frac{2}{t} + \frac{n \, f(\bar{x})^{\frac{1}{n-1}}}{1+t \, f(\bar{x})^{\frac{1}{n-1}}}\] 
for all $t>0$ satisfying $\min_{1 \leq i \leq n} (1+t\lambda_i) > 0$. Integrating this ODE gives 
\[\det P(t) \leq t^2 \, (1+t \, f(\bar{x})^{\frac{1}{n-1}})^n\] 
for all $t > 0$ satisfying $\min_{1 \leq i \leq n} (1+t\lambda_i) > 0$. 

Putting these facts together, we conclude that $\det P(t) = t^2 \, (1+t \, f(\bar{x})^{\frac{1}{n-1}})^n$ for all $t > 0$ satisfying $\min_{1 \leq i \leq n} (1+t\lambda_i) > 0$. Differentiating this identity with respect to $t$, we obtain 
\[\text{\rm tr}(Q(t)) = \frac{2}{t} + \frac{n \, f(\bar{x})^{\frac{1}{n-1}}}{1+t \, f(\bar{x})^{\frac{1}{n-1}}}\] 
for all $t > 0$ satisfying $\min_{1 \leq i \leq n} (1+t\lambda_i) > 0$. Consequently, we must have equality in the arithmetic-harmonic mean equality, and furthermore $\sum_{i=1}^n \lambda_i = n \, f(\bar{x})^{\frac{1}{n-1}}$. Therefore, $\lambda_i = f(\bar{x})^{\frac{1}{n-1}}$ for each $1 \leq i \leq n$. This completes the proof of Lemma \ref{modified.Hessian.of.u.b}. \\

\begin{lemma} 
\label{Hessian.of.u.b}
Assume that $x \in \Omega$. Then $D_\Sigma^2 u(x) = f(x)^{\frac{1}{n-1}} \, g$ and $I\!I(x) = 0$. 
\end{lemma}

\textbf{Proof.} 
By Lemma \ref{modified.Hessian.of.u.b}, $D_\Sigma^2 u(x) - \langle I\!I(x),y \rangle = f(x)^{\frac{1}{n-1}} \, g$ for all $y \in T_x^\perp \Sigma$ satisfying $|\nabla^\Sigma u(x)|^2 + |y|^2 = 1$. Replacing $y$ by $-y$ gives $D_\Sigma^2 u(x) + \langle I\!I(x),y \rangle = f(x)^{\frac{1}{n-1}} \, g$ for all $y \in T_x^\perp \Sigma$ satisfying $|\nabla^\Sigma u(x)|^2 + |y|^2 = 1$. Therefore, $D_\Sigma^2 u(x) = f(x)^{\frac{1}{n-1}} \, g$ and $\langle I\!I(x),y \rangle = 0$ for all $y \in T_x^\perp \Sigma$ satisfying $|\nabla^\Sigma u(x)|^2 + |y|^2 = 1$. From this, the assertion follows easily. \\

\begin{lemma} 
\label{gradient.of.f.b}
Assume that $x \in \Omega$. Then $\nabla^\Sigma f(x) = 0$. 
\end{lemma} 

\textbf{Proof.} 
Let us consider an arbitrary point $x \in \Omega$. Using the definition of $u$, we obtain 
\begin{align*} 
f(x) \, \Delta_\Sigma u(x) 
&= n \, f(x)^{\frac{n}{n-1}} - \sqrt{|\nabla^\Sigma f(x)|^2 + f(x)^2 \, |H(x)|^2} \\ 
&- \langle \nabla^\Sigma f(x),\nabla^\Sigma u(x) \rangle.
\end{align*} 
On the other hand, Lemma \ref{Hessian.of.u.b} implies $\Delta_\Sigma u(x) = n \, f(x)^{\frac{1}{n-1}}$ and $H(x) = 0$. Putting these facts together, we conclude that $\langle \nabla^\Sigma f(x),\nabla^\Sigma u(x) \rangle = -|\nabla^\Sigma f(x)|$. Since $|\nabla^\Sigma u(x)| < 1$, it follows that $\nabla^\Sigma f(x) = 0$. This completes the proof of Lemma \ref{gradient.of.f.b}. \\ 

\begin{lemma}
\label{density.b}
The set $\Omega$ is dense in $\Sigma$. 
\end{lemma}

\textbf{Proof.} 
Suppose that $\Omega$ is not dense in $\Sigma$. Arguing as in Section \ref{proof.of.sobolev.inequality.b}, we obtain 
\[(n+2) \, |B^{n+2}| \, \theta \leq 2 \, |B^2| \int_\Omega f^{\frac{n}{n-1}} < 2 \, |B^2| \int_\Sigma f^{\frac{n}{n-1}} = 2 \, |B^2| \, |B^n| \, \theta.\] 
This contradicts the fact that $(n+2) \, |B^{n+2}| = 2 \, |B^2| \, |B^n|$. This completes the proof of Lemma \ref{density.b}. \\

Since $\Omega$ is a dense subset of $\Sigma$, we conclude that $\nabla^\Sigma f = 0$, $D_\Sigma^2 u = f^{\frac{1}{n-1}} \, g$, and $I\!I = 0$ at each point on $\Sigma$. Since $\Sigma$ is connected, it follows that $f$ is constant. This implies $|\partial \Sigma| = n \, |B^n|^{\frac{1}{n}} \, \theta^{\frac{1}{n}} \, |\Sigma|^{\frac{n-1}{n}}$. 

Note that $u$ is a smooth function on $\Sigma$. Each critical point of $u$ lies in the interior of $\Sigma$ and is nondegenerate with Morse index $0$. In particular, the function $u$ has at most finitely many critical points. 

We next consider the flow on $\Sigma$ generated by the vector field $-\nabla^\Sigma u$. Since the vector field $-\nabla^\Sigma u$ points inward along the boundary $\partial \Sigma$, the flow is defined for all nonnegative times. This gives a one-parameter family of smooth maps $\psi_s: \Sigma \to \Sigma$, where $s \geq 0$. Since $\Sigma$ is connected, standard arguments from Morse theory imply that the function $u$ has exactly one critical point, and $u$ attains its global minimum at that point. It follows that the diameter of $\psi_s(\Sigma)$ converges to $0$ as $s \to \infty$.

Since $D_\Sigma^2 u$ is a constant multiple of the metric, the isoperimetric ratio is unchanged under the flow $\psi_s$. This implies 
\[|\psi_s(\partial \Sigma)| = n \, |B^n|^{\frac{1}{n}} \, \theta^{\frac{1}{n}} \, |\psi_s(\Sigma)|^{\frac{n-1}{n}}\] 
for each $s \geq 0$. If $\theta<1$, this contradicts the Euclidean isoperimetric inequality when $s$ is sufficiently large. Thus, we conclude that $\theta=1$. Consequently, $M$ is isometric to Euclidean space.

Once we know that $M$ is isometric to Euclidean space, the arguments in \cite{Brendle} imply that $\Sigma$ is a flat round ball. This completes the proof of Theorem \ref{rigidity.theorem.b}.


\begin{thebibliography}{99}
\bibitem{Agostiniani-Fogagnolo-Mazzieri}
V.~Agostiniani, M.~Fogagnolo, and L.~Mazzieri, \textit{Sharp geometric inequalities for closed hypersurfaces in manifolds with nonnegative Ricci curvature,} Invent. Math. 222, 1033--1101 (2020)

\bibitem{Allard}
W.~Allard, \textit{On the first variation of a varifold,} Ann. of Math. 95, 417--491 (1972)

\bibitem{Brendle}
S.~Brendle, \textit{The isoperimetric inequality for a minimal submanifold in Euclidean space,} J. Amer. Math. Soc. 34, 595--603 (2021)

\bibitem{Cabre1}
X.~Cabr\'e, \textit{Nondivergent elliptic equations on manifolds with nonnegative curvature,} Comm. Pure Appl. Math. 50, 623--665 (1997)

\bibitem{Cabre2}
X.~Cabr\'e, \textit{Elliptic PDEs in probability and geometry. Symmetry and regularity of solutions,} Discrete Cont. Dyn. Systems A 20, 425--457 (2008)

\bibitem{Cabre-Ros-Oton-Serra}
X.~Cabr\'e, X.~Ros-Oton, and J.~Serra, \textit{Sharp isoperimetric inequalities via the ABP method,} J. Eur. Math. Soc. 18, 2971--2998 (2016)

\bibitem{Carleman} 
T.~Carleman, \textit{Zur Theorie der Minimalfl\"achen,} Math. Z. 9, 154--160 (1921)

\bibitem{Cavalletti-Mondino}
F.~Cavalletti and A.~Mondino, \textit{Sharp and rigid isoperimetric inequalities in metric-measure spaces with lower Ricci curvature bounds,} Invent. Math. 208, 803--849 (2017)

\bibitem{Gilbarg-Trudinger}
D.~Gilbarg and N.~Trudinger, \textit{Elliptic Partial Differential Equations of Second Order,} Springer-Verlag, 2001

\bibitem{Gromov}
M.~Gromov, \textit{Metric structures for Riemannian and non-Riemannian spaces,} Progress in Mathematics vol. 152, Birkh\"auser, Boston, 1999.

\bibitem{Heintze-Karcher}
E.~Heintze and H.~Karcher, \textit{A general comparison theorem with applications to volume estimates for submanifolds,} Ann. Sci. \'Ecole Norm. Sup. 11, 451--470 (1978)

\bibitem{Huisken}
G.~Huisken, \textit{An isoperimetric concept for the mass in general relativity,} Lecture given at the Institute for Advanced Study on March 20, 2009 
\begin{verbatim}
https://www.ias.edu/video/marston-morse-isoperimetric-concept-mass-general-relativity
\end{verbatim}

\bibitem{Klartag}
B.~Klartag, \textit{Needle decompositions in Riemannian geometry,} Mem. Amer. Math. Soc. 249, no. 1180 (2017)

\bibitem{Michael-Simon}
J.H.~Michael and L.M.~Simon, \textit{Sobolev and mean value inequalities on generalized submanifolds of $\mathbb{R}^n$,} Comm. Pure Appl. Math. 26, 361--379 (1973)

\bibitem{Trudinger}
N.~Trudinger, \textit{Isoperimetric inequalities for quermassintegrals,} Ann. Inst. H. Poincar\'e Anal. Non Lin\'eaire 11, 411--425 (1994)

\bibitem{Wang-Zhang}
Y.~Wang and X.W.~Zhang, \textit{Alexandrov-Bakelman-Pucci estimate on Riemannian manifolds,} Adv. Math. 232, 499--512 (2013)

\bibitem{Xia-Zhang}
C.~Xia and X.~Zhang, \textit{ABP estimate and geometric inequalities,} Comm. Anal. Geom. 25, 685--708 (2017)
\end{thebibliography}
\end{document}